\documentclass[11pt]{amsart}
\usepackage{amssymb}
\usepackage{amsfonts}
\usepackage{amscd}
\usepackage{amsmath}
\usepackage{mathrsfs}
\usepackage{curves}
\usepackage{epsfig}
\usepackage[pass]{geometry}
\usepackage{graphicx}
\usepackage{verbatim}
\usepackage{latexsym}
\usepackage{eucal}
\newtheorem{theorem}{Theorem}[section]

\newtheorem{prop}[theorem]{Proposition}

\def \mca {{\mathcal A}}

\def \mch {{\mathcal H}}

\def \mck {{\mathcal K}}

\def \mco {{\mathcal O}}
\def \mcp {{\mathcal Q}}

\def \mcr {{\mathcal R}}

\def \mcv {{\mathcal V}}

\def \mcx {{\mathcal X}}

\def \mbc {{\mathbb C}}

\def \mbn {{\mathbb N}}
\def \mbr {{\mathbb R}}
\def \mbs {{\mathbb S}}

\def \id {\operatorname{Id}}

\def \im {\operatorname{Im}}
\def \re {\operatorname{Re}}
\def \diag{\textrm{Diag}}

\def \beqq {\begin{equation}}
\def \eeqq {\end{equation}}

\def \bpf {\begin{proof}}
\def \epf {\end{proof}}
\def \beq {\begin{equation*}}
\def \eeq {\end{equation*}}

\def \eps {\epsilon}   
\def \la {\lambda}   
\def \La {\Lambda}    
   
\def \lap {\Delta}
\def \p {\partial}
\def \ha {\frac{1}{2}}

\def \rbh {r_{bH}}
\def \rsi {r_{sI}}
\def \xo {\mcx \times_0 \mcx}
\def \ff {\text{ff}}
\def \tilde {\widetilde}

\begin{document}
\title[]{Recovery of  black hole mass from a single quasinormal mode}
\author{Gunther Uhlmann}
\address{Gunther Uhlmann
\newline
\indent Department of Mathematics, University of Washington 
\newline
\indent and Institute for Advanced Study, the Hong Kong University of Science and Technology}
\email{gunther@math.washington.edu}
\author{Yiran Wang}
\address{Yiran Wang
\newline
\indent Department of Mathematics, Emory University}
\email{yiran.wang@emory.edu}
\begin{abstract}
We study the determination of the mass of  a de Sitter-Schwarzschild black hole  from one quasinormal mode. We prove a local uniqueness result with a H\"older type stability estimate. 
\end{abstract}
%\dedicatory{draft on \today}
\date{\today}
 
\maketitle

%%%%%%%%%%%%%%%%%%
\section{Introduction}
The possibility of inferring black hole parameters from quasinormal modes (QNMs) has been explored in the physics literature, see Section 9 of the review paper \cite{BCS}.  For example, for slowly rotating black holes, Detweiler showed by numerical calculation in \cite{Det}  that the wave parameters for the most damped mode are unique functions of the black hole parameters. Later, Echeverria in \cite{Ech} investigated the stability issue. Since the success of gravitational wave interferometers, the topic has gained increasing attention, see e.g.\ \cite{BCW}. One particular motivation for the study is to verify the {\em black hole no hair theorem} for which two QNMs are needed: one QNM is used to recover the black hole parameter and another QNM is used to test the theorem. We refer to \cite[Section 9.7]{BCS} for a review and \cite{Isi} for the state of the art. Despite some convincing evidence, it seems that the theoretical justification is not complete. For example, most of the analysis in the literature is  done for the fundamental modes corresponding to low angular momentum. However, it is generally not known which modes are excited and are extractable from the actual black hole ring down signals, see \cite{BCW, BCS}. In this short note, we aim to  provide a mathematical justification of the recovery of black hole parameters from a single QNM.  

We consider the  model of a non-rotating de Sitter-Schwarzchild black hole $(M, g_{dS})$:
\beqq\label{eq-bh}
\begin{gathered}
M = \mbr_t\times X^\circ, \quad X = (\rbh, \rsi) \times \mbs^2\\
g_{dS} = \alpha^2 dt^2 - \alpha^{-2} dr^2 - r^2 dw^2
\end{gathered}
\eeqq
where $dw^2$ denotes the standard metric on $\mbs^2$ and 
\beqq\label{eq-alpha}
\alpha = (1 - \frac{2m}{r} - \frac{1}{3}\La r^2)^\ha
\eeqq
Here, $m > 0$ is the mass of the black hole and $\La > 0$ is the cosmological constant. They satisfy $0 < 9m^2 \La < 1$. $\rbh, \rsi$ are the two positive roots of $\alpha(r) = 0$ which corresponds to  horizons.  Throughout the note, we assume that $\La$ is known. 
Consider the d'Alembertian on $(M, g_{dS})$: 
\beqq\label{eq-wave}
\square_{M} = \alpha^{-2}(D_t^2 - \alpha^2 r^{-2}D_r(r^2 \alpha^2)D_r - \alpha^2 r^{-2}\lap_{\mbs^2})
\eeqq 
where $D_r = -i \p_r$ and $\lap_{\mbs^2}$ the positive laplacian on $\mbs^2.$ The stationary scattering is governed by the operator 
\beqq\label{eq-wave1}
\lap_X = \alpha^2 r^{-2}D_r(r^2 \alpha^2)D_r + \alpha^2 r^{-2}\lap_{\mbs^2}
\eeqq
see \cite{SaZw}. 
On $L^2(X; \Omega)$ with measure $\Omega = \alpha^{-2}r^2 drdw$, $\lap_X$ is an essentially self-adjoint, non-negative operator, see \cite{MSV1}. 
Consider the resolvent  
\beqq\label{eq-RX}
R_X(\la)  = (\lap_X - \la^2)^{-1}
\eeqq
Here, we use $\la^2$ as the spectral parameter and  take $\im\la \geq 0$ to be the physical plane such that  $R_{X}(\la)$ is  bounded  on $L^2(X; \Omega)$  for $\im\la > > 0$, according to the spectral theorem.  S\'a Barreto and Zworski demonstrated in \cite[Proposition 2.1]{SaZw}  that $R_{X}(\la)$ has a meromorphic continuation as operators from $C_0^\infty(X)$ to $C^\infty(X)$ from $\im \la \geq 0$ to $\mbc$ with poles of finite rank. The poles of $R_X(\la)$ are called resonances. The fact that they are equivalent to the quasinormal modes  defined by using Zerilli's equation (see e.g.\ \cite{ChDe, Cha}) are   discussed in \cite{SaZw}, see also \cite{BoHa, MSV1}.

We denote the set of resonances by $\mcr(m)$ and set $\mco = (0, 1/(3\sqrt{\La}))$. We call $\la$ a trivial resonance if $\la\in \mcr(m)$ for all $m\in \mco.$ For example, it is known that $0$ is a trivial resonance, see \cite{MSV1}. Trivial resonances can not be used to determine black hole parameters. Our main result is 

\begin{theorem}\label{thm-main}
Let $\La > 0$ and $m\in \mco$. For any $\la\in \mcr(m) \backslash i(-\infty, 0]$  not a trivial resonance, there exists $\delta > 0$ (depending on $\la$) such that for any $\tilde m \in \mco$ with $|\tilde m - m| < \delta$,   if $\la \in \mcr(\tilde m)$ then  $m = \tilde m.$ 
Moreover, if $\tilde \la \in \mcr(\tilde m)$ is sufficiently close to $\la$, then 
\beq
|\tilde m - m|  \leq C |\tilde \la - \la|^{1/N}
\eeq
for some $C > 0$ and $N\in \mbn$ depending on $\la$.  
\end{theorem}

We point out that resonances on $i(-\infty, 0]$ are excluded in the theorem. There is a set of resonances on $i(-\infty, 0]$ described by \eqref{eq-mcp} which cannot be treated with our method, although their dependency on $m$ can be found explicitly. We do not investigate it further because such resonances are purely imaginary and they seem to be less relevant in practical cases, see for instance \cite{Isi}. 

\begin{figure}[htbp]
\centering
\includegraphics[scale = 0.61]{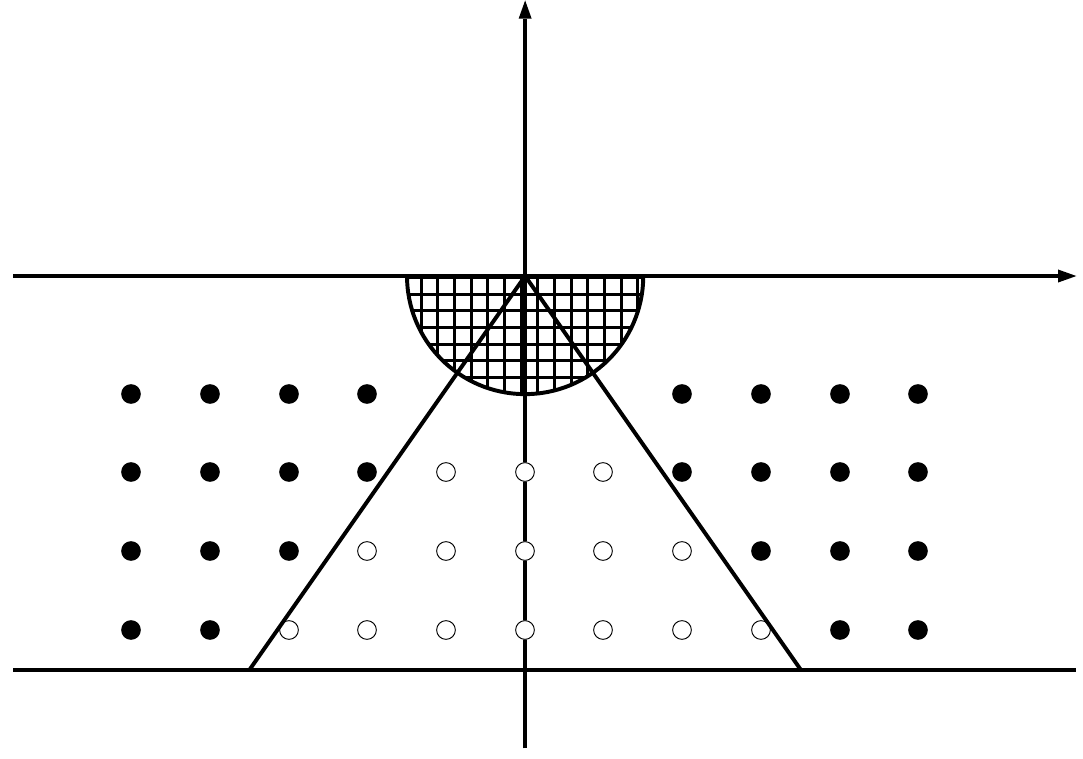}
\caption{Resonances  for de Sitter-Schwarzshild black holes. The black dots are resonances captured by  Theorem in \cite{SaZw}. The hollow dots and resonances in the shaded region are not. }
\label{fig-res}
\end{figure}

We also remark that for recovering black hole parameters, it is common to use only one or a few QNMs. This is very different from the usual inverse spectral/resonance problem for which the whole set $\mcr(m)$ is used to determine the parameters. In fact, there is a large literature on distribution of resonances for large angular momentum. For example, Theorem in \cite{SaZw} states that there exists $K > 0, \theta > 0$ such that for any $C > 0$ there is an injective map $\tilde b$ from the set of pseudo-poles 
\beq
(\pm l \pm \ha - \frac{i}{2}(k + \ha))\frac{(1 - 9 \La m^2)^\ha}{3^{3/2}m}
\eeq
in to $\mcr(m, \La)$ such that all the poles in 
\beq
\Omega_C = \{\la: \im \la > -C, |\la| > K, \im \la > -\theta |\re\la|\}
\eeq
are in the image of $\tilde b$ and for $\tilde b(\mu)\in \Omega_C$, we have $\tilde b(\mu) - \mu \rightarrow 0 \text{ as } |\mu|\rightarrow \infty.$
See Figure \ref{fig-res}. 
By looking at the sequence of resonances for large $l, k,$ one 
recovers $m$. Similar results exist for rotating black holes, see for example \cite{Dya}. 
 Note that resonances sufficiently close to the lattice points  cannot be trivial resonances. %Also, note that the number of resonances in any bounded set of $\mbc$ is finite, so there can only be finitely many points there for which the first alternative applies. 
It is desirable to identify the trivial resonances, if there is any except $0.$ There is some interesting recent result in \cite{HiXi} which shows the convergence of resonances to a set  of $i(-\infty, 0]$ for small masses. The numerical study \cite[Fig. 6(b)]{HiXi} seems to indicate that such points are not trivial resonances.

Our proof of the theorem is based on analytic perturbation argument, by observing that the coefficients of the operator $\lap_X$ are analytic functions in $m$. There are some  resonance perturbation theories, see for instance Agmon \cite{Agm}, Howland \cite{How}, which are developed upon perturbation theory for eigenvalues, see for example \cite{ReSi4}. Here, we use that $\lap_X$ has asymptotically hyperbolic structure near the two horizons to construct a parametrix modulo a trace-class error term, following Mazzeo and Melrose \cite{MaMe}. We then use the  Fredholm determinant and its analyticity in $m$ to finish the proof. This approach has the benefit of not relying on the spherical symmetry of the black hole metric. It is clear from the proof that one can add metric or potential perturbations with suitable decay at the horizons to obtain a similar result to Theorem \ref{thm-main}. 

The note is organized as follows. We begin in Section \ref{sec-examp} with a scattering problem to demonstrate the possibility of recovering parameters from a single resonance. In Section \ref{sec-ahm}, we discuss the asymptotic hyperbolic structure and the analyticity. We construct the resolvent in Section \ref{sec-para} and finish the proof in Section \ref{sec-proof}.

% We also believe that the argument applies to Schwarzschild black holes with $\La = 0$. In that case, the operator $\lap_X$ has an asymptotically hyperbolic structure near the event horizon for which the construction in this note still works, and an asymptotically Euclidean structure near $r = \infty$ for which one can construct a parametrix within Melrose's scattering calculus in \cite{Me}. 

   %%%%%%%%%%%%%%%%%%%%%%%%
 \section{An example: the potential barrier}\label{sec-examp}
In this section, we give an example of a scattering system depending on one parameter, for which a single resonance recovers the parameter. The example was actually used by Chandrasekhar and Detweiler in \cite{ChDe} to illustrate the concept of quasinormal mode. 

\begin{figure}[htbp]
\centering
\vspace{-3.1cm}
\includegraphics[scale = 0.45]{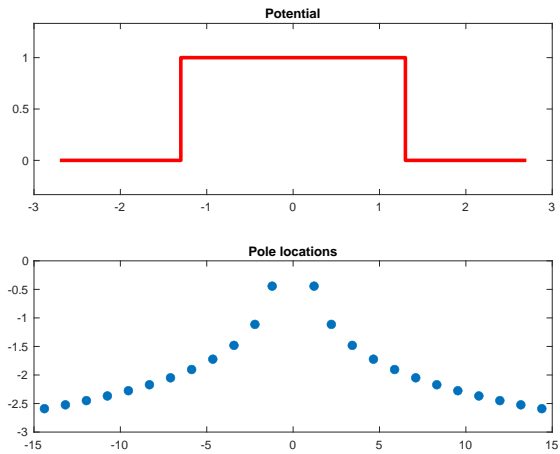} 
\vspace{-3.3cm}
\caption{A rectangular barrier and its resonances.}
\label{fig-pole1}
\end{figure}

Consider  
\beq
u''(x) - V(x) u(x) + \sigma^2 u(x) = 0, \quad x\in \mbr
\eeq
where $\sigma$ is constant and $V$ is the rectangular barrier  
\beq
V(x) = \begin{cases}
1, \quad x\in [-L, L]\\
0, \quad \text{otherwise}
\end{cases}
\eeq
See Figure \ref{fig-pole1}. Note that the potential is characterized by $L$. In this case, the scattering resonances can be defined as the poles of the scattering matrix. It is a standard exercise in scattering theory to find the scattering matrix. Let's look at a wave traveling to the right, hit the potential and gets reflected and transmitted. In this case, the solution looks like
\beq
u_{R}(x) =   \begin{cases} 
  e^{i\sigma x} + r e^{-i\sigma x}, x < -L\\
 t e^{i\sigma x}, x > L
\end{cases}
\eeq
Here, $r$ is   the refection coefficient and $t$ is the transmission coefficient. Similarly, we  can consider a  wave traveling to the left of the form 
\beq
u_{L}(x) =   \begin{cases} 
t' e^{-i\sigma x}, x < -L\\
e^{-i\sigma x} + r' e^{i\sigma x}, x > L
\end{cases}
\eeq
with $r', t'$  the reflection, transmission coefficient respectively. 
%Note that  $u_{R}$ and $u_{L}$ can be expressed in terms of the basic outgoing and incoming waves: 
%\beq
%\begin{gathered}
%u_{in}(x) = e^{i \sigma x}, \quad x <-L \\
%v_{in}(x) = e^{-i\sigma x}, \quad x > L  
%\end{gathered}
%\eeq
%and 
%\beq
%\begin{gathered}
%u_{out}(x) = e^{i \sigma x}, \quad x > L \\
%v_{out}(x) = e^{-i\sigma x}, \quad x < -L  
%\end{gathered}
%\eeq
The scattering matrix is  
% \beq
% \begin{pmatrix}
% u_R\\
% u_L
% \end{pmatrix}
%  = \begin{pmatrix}
% u_{in} \\
%v_{in}
%  \end{pmatrix}
%   + 
%   S \begin{pmatrix}
%u_{out} \\
%v_{out}
%  \end{pmatrix}
% \eeq
%  with 
\beq
S = \begin{pmatrix}
t &r\\
r' &t'
\end{pmatrix}
\eeq
%
%We look for solution $u$ on $[-L, L]$ of the form 
%\beq
%u(x) = Ae^{iq x} + Be^{-iqx} \quad x\in [-L, L]
%\eeq
%where $q^2 = \sigma^2 -1.$ Then we match the solutions at $x = \pm L$ so the solution is $C^1$.   This calculation is a bit long. 
%First we find $r, t. $ At $x = -L,$ we have two equations
%\beq
%\begin{gathered}
%e^{-i\sigma L} + r e^{i\sigma L} = A e^{-iq L} + B e^{iq L}\\
%\sigma e^{-i\sigma L} - r\sigma e^{i\sigma L} = A q e^{-iq L} -  Bq e^{-iqL}
%\end{gathered}
%\eeq
%At $x = L$ we have another two equations
%\beq
%\begin{gathered}
%te^{i\sigma L}   = A e^{iq L} + B e^{-iq L}\\
%t\sigma e^{i\sigma L}  = A q e^{iq L} -  Bq e^{-iqL}
%\end{gathered}
%\eeq
%Eliminating $A, B$, we can solve that 
%\beq
%r = \frac{e^{-2i\sigma L + 2iqL} - e^{-2i\sigma L - 2iqL}}{K}, \quad K = \frac{q + \sigma}{q - \sigma} e^{-2iqL} - \frac{q - \sigma}{q + \sigma} e^{2iqL}
%\eeq
%and 
By matching the solution and its derivatives at $x = L, -L$, we can find the coefficients as 
\beq
\begin{aligned}
&r = r' = \frac{e^{-2i\sigma L + 2iqL} - e^{-2i\sigma L - 2iqL}}{K}, \\
&t  = t' = e^{-2i\sigma L - 2iqL} + r \frac{q + \sigma}{q - \sigma} e^{-2iqL} 
\end{aligned}
\eeq
where 
\beq
K = \frac{q + \sigma}{q - \sigma} e^{-2iqL} - \frac{q - \sigma}{q + \sigma} e^{2iqL}
\eeq  
Thus, the resonances are solutions of $K = 0$ or equivalently 
\beqq\label{eq-sca}
(\frac{q + \sigma}{q - \sigma})^2 = e^{4iqL}
\eeqq
where $q^2 = \sigma^2 - 1.$  Suppose we have $\sigma$ such that $\im q\neq 0$. Then we can take modulus of \eqref{eq-sca} to find $L$ as
\beqq\label{eq-L}
L = -\frac{1}{2\im q} \ln |\frac{q + \sigma}{q - \sigma}|.
\eeqq
This shows that one can recover $L$ from one resonance. 

Now we provide a numerical verification. We compute resonances using a Matlab code from \cite{Code} and  identify $L$ using \eqref{eq-L}. It is important to note that the code from \cite{Code} does not calculate resonances by solving \eqref{eq-L}.  Take $L = 1.3$. The potential and resonances are plotted in Figure \ref{fig-pole1}. The numerical values of the four resonances nearest to the origin with positive real parts are 
 \beq
 \begin{gathered}
\la_1 = 1.2127 - 0.4432i, \quad  \la_2 = 2.2120 - 1.1135i, \\
 \la_3 =   3.4242 - 1.4810i, \quad \la_4 =    4.6501 - 1.7230i  
 \end{gathered}
\eeq
Using any of these resonances in \eqref{eq-L}, we find $L = 1.3$ with a $10^{-4}$ error. 

 %%%%%%%%%%%%%%%%%%%%%%%%
 \section{The asymptoically hyperbolic structure and analyticity}\label{sec-ahm}
It is known that $\lap_X$ in \eqref{eq-wave1} can be essentially viewed as perturbed Laplacians associated with some asymptotically hyperbolic metrics near $\p X$. We follow the presentations in \cite{MSV1}.  
Let $X$ be a compact manifold of dimension $n+1$ with boundary $\p X$. Let $\rho$ be a boundary defining function such that $\rho > 0$ in $X$, $\rho = 0$  at $\p X$, $d\rho \neq 0$  at $\p X.$  
A metric $g$ on $X$ is called conformally compact if $G = \rho^2 g$ is a non-degenerate Riemannian metric on the closure $\overline X$. If in addition $|d\rho|^2_G|_{\p X} = K$ a constant, the metric $g$ is called asymptotically hyperbolic. In this case, the sectional curvature approaches $-K$ along any curve towards $\p X$, see \cite[Lemma (2.5)]{MaMe}. There is a normal form of the metric near $\p X$, see e.g.\ Graham \cite{Gra}. In particular, there is a choice of boundary defining function $x$ such that in a neighborhood $U = [0, \eps)_x\times Y, Y\subset \p X$ of $p\in \p X$, we can use local coordinate $(x, y), y\in Y$ and get  
\beqq\label{eq-metricnormal}
g = \frac{dx^2 + h(x, y, dy)}{x^2} 
\eeqq 

Now we consider $\lap_X$ in \eqref{eq-lapx} on $X$. We define
\beq
\beta = \ha \frac{d \alpha^2}{dr}  = \frac{m}{r^2} - \frac{\La}{3}r
\eeq
We see that $\beta$ is a smooth function of $r$ on $[\rbh, \rsi]$ and analytic in $m \in \mco$.  We set $\beta_{bH} = \beta(\rbh) > 0, \beta_{sI} = \beta(\rsi) < 0.$  Here, we recall that 
 \beqq\label{eq-rpm}
 \begin{gathered}
 \rbh = \im(\sqrt{1 - (3m\sqrt \La)^2} + i3m \sqrt \La)^{1/3} /\sqrt{\La}, \\
 \rsi = \im(-\sqrt{1 - (3m\sqrt \La)^2} + i3m \sqrt \La)^{1/3} /\sqrt{\La},
 \end{gathered}
 \eeqq
 see page 6 of \cite{SaZw}. Thus $\beta_{sI}, \beta_{bH}$ are both analytic functions of $m \in \mco.$  
Now we write \eqref{eq-wave1} as 
 \beqq\label{eq-lapx0} 
 \begin{gathered}
\lap_X %=   \alpha^2 r^{-2}D_r(r^2 \alpha^2)D_r + \alpha^2 r^{-2}\lap_{\mbs^2}\\
 =  \beta r^{-2} \alpha D_\alpha (\beta r^2 \alpha D_\alpha) + \alpha^2 r^{-2}\lap_{\mbs^2}
 \end{gathered}
\eeqq 
For convenience, we denote $\p X = \p X_{sI} \cup \p X_{bH}$ with 
\beq
\p X_{sI} = \{\rsi\}\times \mbs^2, \quad \p X_{bH} = \{\rbh\} \times \mbs^2.
\eeq
Note that $\alpha$ only vanishes at $\p X$. We let $\rho$ be a boundary defining function  defined through 
\beqq\label{eq-rho}
\begin{gathered}
\alpha = 2\rbh\beta_{bH}\rho \text{ near } r = \rbh\\
\text{and } \alpha = 2\rsi\beta_{sI}\rho \text{ near } r = \rsi. 
\end{gathered}
\eeqq
Here, the smooth structure on $X$ is changed. Before, $r - \rbh$ is a smooth boundary defining function near $\p X_{bH}$ but now we think of $(r - \rbh)^\ha$ as a smooth boundary defining function, see \cite[Section 2]{SaZw}. 
By using $\rho$, \eqref{eq-lapx0} becomes
\beqq\label{eq-lapx}
\begin{gathered}
\lap_X = \beta r^{-2} \rho D_\rho (\beta r^2 \rho D_\rho) + 4\rho^2\beta_{bH}^2 \rbh^{2} r^{-2}\lap_{\mbs^2} \text{ near } \p X_{bH}\\
\lap_X = \beta  r^{-2} \rho D_\rho (\beta r^2 \rho D_\rho) + 4\rho^2\beta_{sI}^2 \rsi^{2} r^{-2}\lap_{\mbs^2} \text{ near } \p X_{sI}
\end{gathered}
\eeqq
Let $g_{bH}$ be the metric defined in a neighborhood of $\p X_{bH}$ given by
 \beqq\label{eq-gbh}
 \begin{gathered}
 g_{bH} = \frac{d\rho^2}{\beta^2 \rho^2} + \frac{r^2}{(2 \beta_{bH} \rbh)^{2}} \frac{dw^2}{\rho^2} 
 \end{gathered}
 \eeqq
 and let $g_{sI}$ be the metric defined in a neighborhood of $\p X_{sI}$ given by
 \beqq\label{eq-gsi}
 \begin{gathered}
 g_{sI} = \frac{d\rho^2}{\beta^2 \rho^2} + \frac{r^2}{(2 \beta_{sI} \rsi)^{2}} \frac{dw^2}{\rho^2}
 \end{gathered}
 \eeqq
 These can be viewed as metric perturbations of the  hyperbolic metrics 
 \beq
 \begin{gathered}
 g_{bH, 0} = \frac{4dz^2}{\beta_{bH}^2 (1 - |z|^2)} \quad g_{sI, 0} = \frac{4dz^2}{\beta_{sI}^2 (1 - |z|^2)}
 \end{gathered}
 \eeq 
 on $\mathbb{B}^{3} = \{z\in \mbr^3: |z| \leq 1\}$ with constant negative sectional curvature $-\beta^2_{bH}$  and $-\beta^2_{sI}$ respectively.  Here, $(1 - |z|^2)^\ha$ is the boundary defining function.  Also, $g_{bH}, g_{sI}$ are even asymptotically hyperbolic metrics as defined in Guillarmou \cite{Gui}. 
 
After some calculation see \cite[Proposition 8.1]{MSV1}, we conclude that there are two smooth functions $W_{bH}, W_{sI}$  such that 
 \beqq\label{eq-lapx2}
 \begin{gathered}
\alpha\lap_X\alpha^{-1}  =  \rho \lap_X \rho^{-1} = \lap_{g_{bH}} + \rho^2 W_{bH} - \beta_{bH}^2, \text{ near }\p X_{bH}\\
\alpha\lap_X\alpha^{-1} =   \rho \lap_X \rho^{-1} = \lap_{g_{sI}} + \rho^2 W_{sI} - \beta_{sI}^2, \text{ near } \p X_{sI}
 \end{gathered}
 \eeqq
This shows the asymptotically hyperbolic structure of $\alpha\lap_X\alpha^{-1}$ near $\p X.$ 

Consider the dependency of the operator $\alpha\lap_X \alpha^{-1}$ on $m$. Note that the manifold $X$   varies when varying $m$. We change the notation from $X$ to $X(m)$. The dependency can be fixed by transforming $X(m)$ to a fixed reference manifold. Let $\mcx = (1, 2) \times \mbs^2$ and let 
\beq
\Psi: \mcx \rightarrow X(m) 
\eeq
 be a diffeomorphism defined by $(s, w) = \Psi(s, w) = (\psi(s), w)$ with 
 \beq
\psi(s) = (s - 1) \rsi + (2- s)\rbh
 \eeq
 Note that $\Psi$ extends smoothly to $\overline \mcx \rightarrow \overline{X(m)}$. 
 Since $\rbh, \rsi$ are analytic functions of $m\in \mco$,  $\psi$ and $\Psi$ are also analytic in $m \in \mco.$ 
Now, the pull back $\rho^* = \psi^*(\rho)$ is a family of smooth boundary defining functions for $\p \mcx$. To see their dependency on $m$, we write \eqref{eq-alpha} as 
\beq
\alpha = \sqrt{\La/3} (r - \rbh)^\ha(r - \rsi)^\ha(r  - r_0)^\ha
\eeq
where $\rbh, \rsi$ are two positive roots of $\alpha = 0$ and $r_0$ is the third negative root. Using \eqref{eq-rho} and on $\mcx$ near $\p \mcx_{bH}$, we have 
\beqq
\begin{aligned}
\rho^* &=  \frac{\sqrt{\La/3} (r - \rbh)^\ha(r - \rsi)^\ha(r  - r_0)^\ha}{2\beta_{bH}\rbh}\\
 &= \frac{\sqrt{\La/3}(\rsi - \rbh) (s- 1)^\ha(2 - s)^\ha((s - 1)\rsi + (2 - s)\rbh - r_0)^\ha}{2\beta_{bH}\rbh}\\
 &= (s - 1)^\ha A_{bH}(s, m)
\end{aligned}
\eeqq 
where $A_{bH}(s, m)$ is defined through the last two lines. It is clear that $A_{bH}$ is smooth in $s$. Since $(s - 1)\rsi + (2 - s)\rbh - r_0 > 0$ for $s\in [1, 2]$, we see that $A_{bH}$, hence $\rho^*$, is analytic in $m\in \mco.$ Near $\p \mcx_{sI}$, we have a similar form 
\beqq\label{eq-rho1}
\begin{gathered}
\rho^*   = (2-s)^\ha A_{sI}(s, m).
\end{gathered}
\eeqq
From \eqref{eq-gbh}, \eqref{eq-gsi}, we see that the pull-back of the metrics are   
 \beqq\label{eq-g1}
 \begin{gathered}
 \Psi^*g_{bH} = \frac{(d\rho^*)^2}{(\rsi - \rbh)^2(\beta^*)^2 (\rho^*)^2} + \frac{[(s - 1)\rsi + (2-  s)\rbh]^2}{(2 \beta_{bH} \rbh)^{2}} \frac{dw^2}{(\rho^*)^2} \\ 
\Psi^* g_{sI} =\frac{(d\rho^*)^2}{(\rsi - \rbh)^2(\beta^*)^2 (\rho^*)^2} + \frac{[(s - 1)\rsi + (2-  s)\rbh]^2}{(2 \beta_{sI} \rsi)^{2}} \frac{dw^2}{(\rho^*)^2} 
 \end{gathered}
 \eeqq
 near $\p X_{bH}, \p X_{sI}$ respectively. Here, $\beta^* = \psi^*(\beta).$

 Let $\mcv_0(\mcx)$ be the Lie algebra of smooth vector fields on $\mcx$ vanishing  at $\p \mcx$.  In local coordinate $(x, y)$ near $\p \mcx$ with $x$ being the boundary defining function, $\mcv_0(\mcx)$ is generated by $x\p_x, x\p_y$. The space of zero-differential operators of order $m$ on $\mcx$, denoted by $\text{Diff}_0^m(\mcx)$ is generated by $m$ fold products of vector fields in $\mcv_0(\mcx)$.  From the  analyticity of the diffeomorphism $\Psi$ on $m$, we see that the pull-back of $\alpha\lap_{X}\alpha^{-1}$ to $\mcx$ is a differential operator on $\mcx$ with coefficients analytic in $m.$ To see it belongs to  $\text{Diff}_0^2(\mcx)$ with coefficients analytic in $m \in \mco$,  we consider  $\alpha \lap_X \alpha^{-1}$  near $\p \mcx$ for example near $s = 1$. We change the boundary defining function from $\rho^*$ to $\gamma = (s - 1)^\ha$ so $\rho^* = \gamma A_{bH}(\gamma^2, m)$. Then the metric in \eqref{eq-g1} becomes
\beq
\Psi^*g_{bH} =   \frac{(1 + \p_\gamma A_{bH})(d \gamma)^2}{(\rsi - \rbh)^2(\beta^*)^2 A^2_{bH} \gamma^2} + \frac{[(s - 1)\rsi + (2-  s)\rbh]^2}{(2 \beta_{bH} \rbh)^{2} A_{bH}^2} \frac{dw^2}{\gamma^2}
\eeq
which is a family of Riemannian metrics on $\mcx$ analytic in $m.$ 

 %%%%%%%%%%%%%%%%%%%
 \section{The resolvent construction}\label{sec-para}
We obtain an approximation of $R_X(\la)$ in \eqref{eq-RX} following Mazzeo-Melrose \cite{MaMe}. In fact, we will find the resolvent of $\alpha \lap_X \alpha^{-1}$ on $\mcx$.   We will be using operators acting on half densities on $\mcx$. For convenience, we introduce an auxiliary  Riemannian metric $g_\mcx$ on $\mcx$ which equals $g_{bH}, g_{sI}$ near $\p \mcx_{bH}, \p \mcx_{sI}$ respectively. Such a metric can be obtained by gluing $g_{bH}, g_{sI}$ near $\p \mcx$ and some Riemannian metric in the interior of $\mcx$. The choice is clearly not unique and its dependency on $m$ is not important.  We use $g_\mcx$ to trivialize the (zero) one-density bundle $\Omega_0$, that is we take $\Omega_0$ to be the volume form $|dg_\mcx|$. The half-density bundle is $\Omega_0^\ha.$ 
 Let $x$ be a boundary defining function  such that in local coordinates $(x, y), x \geq 0, y\in \mbs^2$ near $\p \mcx$, $g_\mcx$ is expressed in form of  \eqref{eq-metricnormal}. In this coordinate, 
\beq
 \Omega^\ha_0 = H(x, y)|\frac{dx}{x}\frac{dy}{x}|^\ha
\eeq
for some smooth function $H$. 
Now we consider $\alpha \lap_X \alpha^{-1}$ acting on smooth sections $C^\infty(\mcx; \Omega^\ha_0)$  in the following way
\beq
\alpha \lap_X \alpha^{-1} (u \Omega^\ha_0) = (\alpha \lap_X \alpha^{-1}u) \Omega^\ha_0
\eeq
The resolvent $R_\alpha(\la) = (\alpha \lap_X \alpha^{-1} - \la^2)^{-1}$ is  acting on $\Omega^\ha_0$ in the same way.

\begin{figure}[htbp]
\centering
\includegraphics[scale = 0.57]{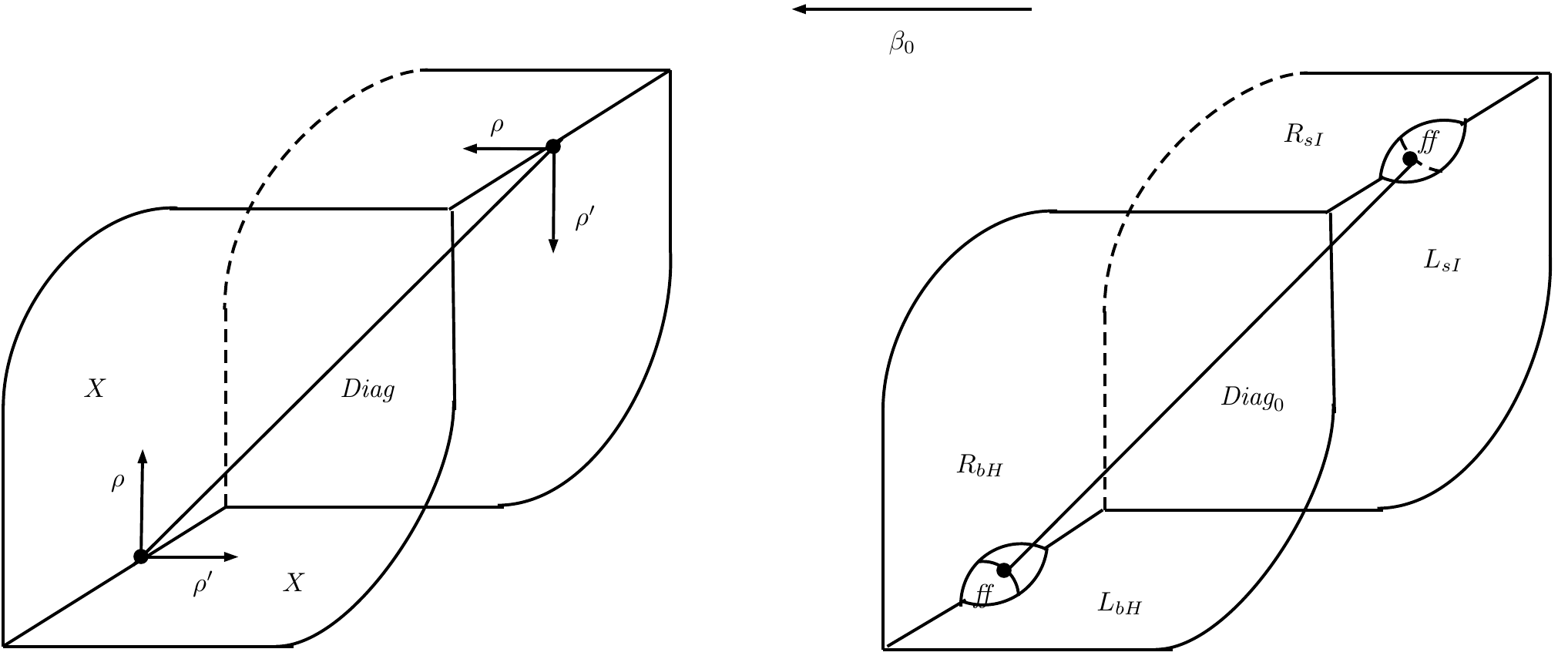}
\caption{The 0-blown up space. The blown-up at the two components of $\p \diag_0$ are shown.}
\label{fig-0blow}
\end{figure}

The parametrix is constructed on the $0$-blown-up space of $\mcx\times \mcx$ as in \cite{MaMe}. Let $\diag = \{(z, z)\in \mcx\times \mcx\}$ be the diagonal of $\mcx\times \mcx$. Let  $\p \diag = \diag \cap(\p \mcx\times \p \mcx)$ which has two (disjoint) connected components.  
As a set, the $0$-blown-up space  is 
\beq
 \mcx \times_0 \mcx = (\mcx\times \mcx)\backslash \p\diag \sqcup S_{++}(\p \diag),
\eeq
where $S_{++}(\p\diag)$ denotes the inward pointing spherical bundle of $T_{\p\diag}^*(\mcx \times \mcx)$. Let 
\begin{gather}
\beta_0: \mcx \times_0 \mcx \rightarrow \mcx\times \mcx \label{0blowdown}
\end{gather}
be the blow-down map. Then $\mcx \times_0 \mcx$ is equipped with a topology and smooth structure  of a manifold with corners for which $\beta_0$ is smooth. The manifold $\xo$  has the following boundary hyper-surfaces: the left and right faces $L=\overline{\beta_0^{-1}(\p \mcx \times \mcx)},$  $R=\overline{\beta_0^{-1}(\mcx \times \p \mcx)},$  and the front face $\ff= \overline{\beta_0^{-1}(\p\diag)}$. Since $\p \mcx = \p \mcx_{bH}\cup \p \mcx_{sI}$ where the asymptotic behavior of the resolvent is different at each connected component, it is convenient to introduce 
\beq
\begin{gathered}
L_{bH} =\overline{\beta_0^{-1}(\p \mcx_{bH} \times \mcx)}, \quad L_{sI} =\overline{\beta_0^{-1}(\p \mcx_{sI} \times \mcx)},\\
R_{bH} =\overline{\beta_0^{-1}(\mcx \times \p \mcx_{bH})}, \quad R_{sI} =\overline{\beta_0^{-1}(\mcx \times \p \mcx_{sI})}
\end{gathered}
\eeq
so $L = L_{bH}\cup L_{sI}$, $R = R_{bH}\cup R_{sI}$, see Figure \ref{fig-0blow}. 
The lifted diagonal is denoted by $\diag_0 = \overline{\beta_0^{-1}(\diag\setminus \p \diag)}$.  $\xo$ has co-dimension two corners at   the intersection of  two of the boundary faces $L, R, \ff$ and co-dimension three corners given by the intersection of all the three faces. See Figure \ref{fig-0blow}.

Now we introduce spaces of operator  on $\xo$. First, let $\mck_0^m(\mcx) \subset \mathscr{D}'(\xo; \Omega_0^\ha)$ be the space of conormal distributions of the bundle $\Omega_0^\ha$ to $\diag_0$ and vanishing to infinite orders at $L, R$. Here, it is understood that $\Omega_0^\ha$ denotes the half-density bundle lifted from the one on $\mcx\times \mcx$ by $\beta_0.$ The corresponding class of pseudo-differential operators is denoted by $\Psi_0^m(\mcx, \Omega_0^\ha)$. Next, let $\mathscr{V}_b$ be the space of smooth vector fields on $\xo$  which are tangent to each of the boundary faces $L, R, \ff$. Let $\rho_{\bullet}, \bullet = L_{bH}, L_{sI}, R_{bH}, R_{sI}, \ff$ be boundary defining functions. We set 
\beqq\label{eq-spaceA}
\begin{gathered}
\mca^{a, b, c, d}(\xo) = \{u \in \mathscr{D}'(\xo): V_1\cdots V_k u\in \\
\rho_{L_{bH}}^a\rho_{R_{bH}}^b \rho_{L_{sI}}^c \rho_{R_{sI}}^d L^\infty(\xo), V_i \in \mcv_b, i = 1, 2, \cdots, k, \forall k \geq 0\} 
 \end{gathered}
\eeqq
Then define
\beq
\mck_0^{-\infty, a, b, c, d}(\xo) = \mca^{a, b, c, d}(\xo)\otimes C^\infty(\xo; \Omega_0^\ha).
\eeq
Finally, define
\beq
\mck^{m, a, b, c, d}_0(\mcx) = \mck_0^{-\infty, a, b, c, d}(\xo) + \mck_0^m(\mcx).
\eeq
Then we let $\Psi_0^{m, a, b, c, d}(\mcx)$ be the space of operators on $\mcx$ whose Schwartz kernel when lifted to $\xo$ belongs to $\mck^{m, a, b, c, d}_0(\mcx)$.  

We have the following result. 
 \begin{prop}
  %Consider the operator $\alpha \lap_X \alpha^{-1}$. %and the resolvent $R_\alpha(\la) = (\alpha \lap_X \alpha^{-1} - \la^2)^{-1}$ as a bounded operator on $L^2(\mcx)$ for $\im\la >>0$. 
  There is a family of operators $M(\la, m) \in \Psi_0^{-2, a, a, b, b}(\mcx)$  with 
  \beqq
  a = 1 + \frac{\la}{\beta_{bH}}i, \quad b = 1 + \frac{\la}{|\beta_{sI}|} i,
  \eeqq
 analytic in $m\in \mco$ and holomorphic in $\la\in \mbc\backslash \mcp$ with 
  \beqq\label{eq-mcp}
  \mcp \doteq \frac{-i}{\beta_{bH}} \mbn \cup \frac{-i}{|\beta_{sI}|}\mbn \cup \{0\}
  \eeqq
  % and residues of finite rank 
  such that 
  \beqq\label{eq-para}
 (\alpha \lap_X \alpha^{-1} - \la^2) M(\la, m) = \id + E(\la, m)
  \eeqq
  Here, $E(\la, m) \in \rho_{\ff}^\infty \Psi_0^{-\infty, \infty, a, \infty, b}(\mcx)$ is trace class on $x^l L^2(\mcx)$ for any $l.$ Moreover, its Schwarz kernel is holomorphic in $\la\in \mbc\backslash \mcp$, and analytic in $m\in \mco$. 
 \end{prop}
 \bpf
For fixed $m$, the construction of $M(\la, m)$ and $E(\la, m)$ and their holomorphy in $\la$ is essentially contained in Proposition (7.4) of \cite{MaMe}, which applies to the laplacian of asymptotically hyperbolic metrics. As argued in Proposition 2.2 of \cite{SaZw}, the result applies to $\alpha \lap_X\alpha^{-1}-\la^2$ as the normal operator is elliptic.  Because we argued in Section \ref{sec-ahm} that $\alpha \lap_X\alpha^{-1} \in \text{Diff}_0^2(\mcx)$ with coefficients analytic in $m\in \mco$, the construction in \cite{MaMe} produces $M(\la, m), E(\la, m)$ analytic in $m.$ The set $\mcp$ comes from the poles of the resolvent of 
\beq
\beta_{bH}^{-2} \lap_{0} -\la^2 
\text{ and } \beta_{sI}^{-2} \lap_{0} -\la^2
\eeq
see Lemma (6.15) of \cite{MaMe}. Here, $\lap_0$ denotes the positive laplacian of the standard hyperbolic metric. By rescaling the operator, we find the set $\mcp.$

To find $a, b$, it is convenient to work on $X$. It suffices to consider the operators near $\p X_{bH}, \p X_{sI}$ respectively. Write $g_{bH} = \rho^{-2}h_{bH}$. From \cite[Theorem (7.1)]{MaMe}, see also \cite[Theorem 1.1]{Gui}, we know that the resolvent  of 
\beqq\label{eq-lap3} 
|d\rho|_{h_{bH}}^{-2} \lap_{g_{bH}} + \zeta(\zeta -2)
\eeqq
belongs to $\Psi_0^{-2, \zeta, \zeta}(\mcx)$. Here, we followed \cite{MaMe} and used a different spectral parameter $\zeta$. Near $\p \mcx_{bH}$, $-|d\rho|_{h_{bH}}^2$ approaches $-\beta_{bH}^2$.     Comparing \eqref{eq-lap3} with \eqref{eq-lapx}, we get 
\beq
\beta_{bH}^2 \zeta(\zeta -2) = -\la^2 - \beta_{bH}^2
\eeq  
which gives $\zeta = 1 + i \la/\beta_{bH}$. This gives $a$, and $b$ can be found in the same way near $\p X_{sI}.$

To see the trace class property, we recall a result \cite[Lemma (5.24)]{MaMe}. The push forward of the space $\rho_{\ff}^\infty\mck_0^{-\infty, a, b, c, d}(\xo)$ is  
\beq
\mca_0^{a, b, c, d}(\mcx\times \mcx; \Omega_0^\ha\otimes\Omega_0^\ha) = \bigcap_{p} (|y - y'|^2 + \rho^2 + (\rho')^2)^p \mca^{a, b, c, d}(\mcx\times \mcx)
\eeq
with the latter defined similarly to \eqref{eq-spaceA}.  Let $K_E(z, z')$ be the Schwarz kernel of $E$. As $E \in \rho_{\ff}^\infty \Psi_0^{-\infty, \infty, a, \infty, b}(\mcx)$, we conclude that for $N\in \mbn$ we can write
\beq
K_E(z, z, m)  = \rho^N  F_N(z, z, m)
\eeq
where $F_N \in C^\infty(\overline \mcx)$ is analytic in $m\in \mco$. Near $\p \mcx_{bH}$, using \eqref{eq-gbh} and $z = (x, y)$ as local coordinate, we see that $dg(z) = H(x, y)x^{-3}|d\rho dy|$.  
We get a similar expression near $\p \mcx_{sI}$. Then we see that the integral  $\int_\mcx |K_E(z, z)|dg_\mcx(z)$ is finite so $E$ is of trace class. 
 \epf

%%%%%%%%%%%%%%%%%
\section{Proof of Theorem \ref{thm-main}}\label{sec-proof}
We apply the resolvent to \eqref{eq-para} to get 
 \beqq
 M(\la, m) = R_\alpha(\la) (\id + E(\la, m))
 \eeqq
 Since $E(\la, m)$ is compact on $x^l L^2(\mcx)$, using analytic Fredholm theorem, we see that for any $m\in \mco$, $(\id + E(\la, m))^{-1}$ is a family of bounded operators, meromorphic in $\la \in \mbc\backslash \mcp$. The poles (at least away from $\mcp$) are the resonances. In fact, the resolvent $R_\alpha(\la)$ is also meromorphic at $\mcp$, as clarified in \cite{Gui}.

 Now we use the determinant of $\id + E$ to analyze the poles.  We recall that if $A$ is a trace class operator on a Hilbert space $\mch$ with eigenvalues $\la_k, k = 1, 2, \cdots$ with $|\la_1|\geq |\la_2|\geq \cdots \geq 0$. Then the Fredholm determinant $\det(\id + A) = \Pi_{k = 1}^\infty (1 + \la_k)$. See \cite[Appendix B]{DyZw}. Also, $\id + A$ is invertible if and only if $\det(\id + A)$ is non-zero, see \cite[Proposition B.28]{DyZw}. Therefore, the resonances of $R_\alpha(\la)$ is contained in the zero set of 
\beq
K(\la, m) = \det (\id + E(\la, m))
\eeq
Using the argument in the end of \cite[Section B.5]{DyZw}, we conclude that $K(\la, m)$ is a function holomorphic in $\la \in \mbc\backslash \mcp$, and analytic in $m\in \mco.$ 

Now we   suppose $\la_0$ is a resonance so $K(\la_0, m) = 0$. By the analyticity in $m$, either $K(\la_0, m)$ is identically zero for all $m$ which means $\la_0$ is a resonance for all $m\in \mco$, or $m$ is the only (discrete) zero locally.  This proves the first claim of Theorem \ref{thm-main}.

For the stability, we write $K(\la, m) = (m - m_0)^N f(\la, m)$ for some $N\geq 0$ and $f$ analytic in $m$ with $f(\la_0, m_0)\neq 0$. Now we set $t = (m - m_0)^N$ and get 
\beq
K(\la, t) = K(\la, m) = t f(\la, m_0 + t^{1/N}). 
\eeq
We see that 
\beq
\p_t K(\la, t)|_{t = 0} = f(\la, m_0)\neq 0
\eeq
Using the implicit function theorem, %(see e.g.\ \cite[Theorem 2-12]{Spi}), 
we get that $t = g(\la)$ is differentiable in a neighborhood of $\la_0$. Thus,  $|t|\leq C |\la -\tilde \la|$ 
which implies 
\beq
|m - \tilde m|\leq C|\la - \tilde \la|^{1/N}
\eeq
 This completes the proof of the Theorem \ref{thm-main}.

\section*{Acknowledgment}
The authors thank Peter Hintz for his  thorough  reading of a previous version of the manuscript and for making very useful comments.  GU was partly supported by NSF, a
Simons Fellowship, a Walker Family Endowed Professorship at UW and a Si-Yuan Professorship at IAS, HKUST.

%%===============================REFERENCE==========================================%


\begin{thebibliography}{99}
\bibitem{Agm} S. Agmon. {\em A perturbation theory of resonances.} Communications on Pure and Applied Mathematics 51.11‐12 (1998): 1255-1309.

\bibitem{BCS} E. Berti, V. Cardoso, A. Starinets. {\em Quasinormal modes of black holes and black branes.} Classical and Quantum Gravity 26.16 (2009): 163001.

\bibitem{BCW} E. Berti, V. Cardoso, C. Will. {\em Gravitational-wave spectroscopy of massive black holes with the space interferometer LISA.} Physical Review D 73.6 (2006): 064030.

\bibitem{BoHa} J.-F. Bony, D. H\"afner. {\em Decay and non-decay of the local energy for the wave equation on the de Sitter–Schwarzschild metric.} Communications in Mathematical Physics 282.3 (2008): 697-719.

\bibitem{Cha} S. Chandrasekhar. {\em The mathematical theory of black holes.} The International Series of Monographs on Physics,  Volume 69, Clarendon Press, Oxford, UK (1983). 

\bibitem{ChDe} S. Chandrasekhar, S.  Detweiler. {\em The quasi-normal modes of the Schwarzschild black hole.} Proceedings of the Royal Society of London. A. Mathematical and Physical Sciences 344.1639 (1975): 441-452.

\bibitem{Det} S. Detweiler. {\em Black holes and gravitational waves. III -- The resonant frequencies of rotating holes.}  The Astrophysical Journal 239 (1980): 292-295.

\bibitem{Dya} S. Dyatlov. {\em Asymptotic distribution of quasi-normal modes for Kerr–de Sitter black holes.}  Annales Henri Poincaré. Vol. 13. No. 5. SP Birkh\"auser Verlag Basel, 2012.

\bibitem{DyZw} S. Dyatlov, M. Zworski. {\em Mathematical theory of scattering resonances.} Vol. 200. American Mathematical Soc., 2019.

\bibitem{Ech} F. Echeverria. {\em Gravitational-wave measurements of the mass and angular momentum of a black hole.} Physical Review D 40.10 (1989): 3194.

\bibitem{Gra} R. Graham. {\em Volume and area renormalizations for conformally compact Einstein metrics.} Proceedings of the 19th Winter School. Circolo Matematico di Palermo, 2000.

\bibitem{Gui} C. Guillarmou. {\em Meromorphic properties of the resolvent on asymptotically hyperbolic manifolds.} Duke Mathematical Journal 129.1 (2005): 1-37.

\bibitem{HiXi} P. Hintz, Y. Xie. {\em Quasinormal modes of small Schwarzschild–de Sitter black holes.} Journal of Mathematical Physics 63.1 (2022): 011509.


\bibitem{How} J. Howland. {\em Puiseux series for resonances at an embedded eigenvalue.} Pacific Journal of Mathematics 55.1 (1974): 157-176.

\bibitem{Isi} M. Isi, M. Giesler, W. Farr, M. Scheel, S.  Teukolsky (2019). {\em Testing the no-hair theorem with GW150914.} Physical Review Letters, 123(11), 111102.
 
%\bibitem{JoSa} M. Joshi, A. S\'a Barreto. {\em Inverse scattering on asymptotically hyperbolic manifolds.} Acta Mathematica 184.1 (2000): 41-86.

\bibitem{MaMe} R. Mazzeo, R. Melrose. {\em Meromorphic extension of the resolvent on complete spaces with asymptotically constant negative curvature.} Journal of Functional Analysis 75.2 (1987): 260-310.

%\bibitem{Me} R. Melrose. {\em Spectral and scattering theory for the Laplacian on asymptotically Euclidean spaces.} Lecture Notes in Pure and Applied Mathematics (1994): 85-85.

\bibitem{MSV1} R. Melrose, A. S\'a Barreto, A. Vasy. {\em Analytic continuation and semiclassical resolvent estimates on asymptotically hyperbolic spaces.} Communications in Partial Differential Equations 39.3 (2014): 452-511.

%\bibitem{MSV2} R. Melrose, A. S\'a Barreto, A. Vasy. {\em Asymptotics of solutions of the wave equation on de Sitter-Schwarzschild space.} Communications in Partial Differential Equations 39.3 (2014): 512-529.

\bibitem{ReSi4} M. Reed, B. Simon. {\em Methods of modern mathematical physics IV: Analysis of operators.} New York: Academic. (1978).
 
\bibitem{SaZw} A. S\'a Barreto, M. Zworski. {\em Distribution of resonances for spherical black holes.} Mathematical Research Letters 4.1 (1997): 103-121.

%\bibitem{Spi} M. Spivak. {\em Calculus on manifolds: a modern approach to classical theorems of advanced calculus.} CRC press, 2018.

\bibitem{Code} https://www.cs.cornell.edu/bindel/cims/resonant1d/

\end{thebibliography}
\end{document}